\documentstyle[12pt, fullpage, amsfonts, euscript, psfig]{amsart}
\newtheorem{lemma}{Lemma}
\newtheorem{corollary}{Corollary}

\newenvironment{pflike}[1]{\noindent{\bf #1}}{\qed\medskip}
\newenvironment{proof}{\begin{pflike}{Proof:}}{\end{pflike}}

\newcommand{\story}[1]{{\smallskip\leftskip=26pt\rightskip=\leftskip\noindent\smaller\smaller\baselineskip=12pt#1\par\smallskip}}
\renewcommand{\put}[3]{\setbox0=\hbox{#1}\dimen0=\wd0\hskip#2\raise#3\box0\hskip-#2\hskip-\dimen0}

\newcommand{\oecases}[3]{\begin{cases}#1,&\hbox{if $h$ is odd,}\\#2,&\hbox{if $h$ is even#3}\end{cases}}

\newcommand{\floor}[1]{\lfloor#1\rfloor}
\newcommand{\bfloor}[1]{\big\lfloor#1\big\rfloor}
\newcommand{\bbfloor}[1]{\bigg\lfloor#1\bigg\rfloor}
\newcommand{\ceil}[1]{\lceil#1\rceil}

\newcommand{\rmif}{\hbox{if\ }}

\renewcommand{\cal}{\EuScript}
\newcommand{\A}{{\cal A}}
\newcommand{\N}{{\mathbb N}}
\newcommand{\R}{{\mathbb R}}

\newcommand{\tms}{\ifmmode\times\else$\times$\fi}
\renewcommand{\and}{\hbox{ and }}
\renewcommand{\sp}{semi\-pe\-rim\-e\-ter}
\newcommand{\sps}{\sp s}
\newcommand{\asq}{al\-most\discretionary{-}{}{-}square}
\newcommand{\asqs}{\asq s}

\begin{document}
\title{Farmer Ted Goes Natural}
\author{Greg Martin}
\address{Department of Mathematics\\University of Toronto\\Canada M5S 3G3}
\email{gerg@@math.toronto.edu}
%\subjclass{11B83}
\maketitle

\section{Setting the Stage}\label{introsec}\noindent
We've all been given a problem in a calculus class remarkably similar
to the following one:

\story{Farmer Ted is building a chicken coop. He decides he can spare
190 square feet of his land for the coop, which will be built in the
shape of a rectangle. Being a practical man, Farmer Ted wants to spend
as little as possible on the chicken wire for the fence. What dimensions
should he make the chicken coop?}

\noindent By solving a simple optimization problem, we learn that
Farmer Ted should make his chicken coop a square with side lengths
$\sqrt{190}$ feet. And that, according to the solution manual, is
that.

But the calculus books don't tell the rest of the story:

\story{So Farmer Ted went over to Builders Square and told the
salesman, ``I'd like $4\sqrt{190}$ feet of chicken wire, please.'' The
salesman, however, replied that he could sell one foot or two feet or
a hundred feet of chicken wire, but what the heck was $4\sqrt{190}$
feet of chicken wire? Farmer Ted was taken aback, explaining heatedly
that his family had been buying as little chicken wire as possible for
generations, and he really wanted $4\sqrt{190}$ feet of chicken wire
measured off for him immediately! But the salesman, fearing more
irrational behavior from Farmer Ted, told him, ``I don't want to hear
about your roots. We do business in a natural way here, and if you
don't like it you can leave the whole store.'' Well, Farmer Ted didn't
feel that this treatment was commensurate with his request, but he
left Builders Square to rethink his coop from square one.}

At first, Farmer Ted thought his best bet would be to make a
$10'\tms19'$ chicken coop, necessitating the purchase of 58 feet of
chicken wire---certainly this was better than 86 feet of chicken wire
for a $5'\tms38'$ coop, say. But then he realized that he could be
more cost-effective by not using all of the 190 square feet of land he
had reserved for the coop. For instance, he could construct an
$11'\tms17'$ coop (187 square feet) with only 56 feet of chicken
wire; this would give him about 3.34 square feet of coop space per
foot of chicken wire purchased, as opposed to only 3.28 square feet per
chicken-wire-foot for the $10'\tms19'$ coop. Naturally, the
parsimonious farmer wondered: could he do even better?

\section{Posing the Problem}\label{problemsec}\noindent
Jon Grantham posed the following problem at the 1998 SouthEast
Regional Meeting On Numbers in Greensboro, North Carolina: given a
positive integer $N$, find the dimensions of the rectangle with
integer side lengths and area at most $N$ whose area-to-perimeter
ratio is largest among all such rectangles. In the story above, Farmer
Ted is trying to solve this problem for $N=190$.

Let's introduce some notation so we can formulate Grantham's problem
more precisely. For a positive integer $n$, let $s(n)$ denote the
least possible \sp\ (length plus width) of a rectangle with integer
side lengths and area $n$. (Since the area-to-\sp\ ratio of a
rectangle is always twice the area-to-perimeter ratio, it doesn't
really change the problem if we consider \sps\ instead of perimeters;
this will eliminate annoying factors of 2 in many of our formulas.) In
other (and fewer) words,
\begin{equation*}
s(n) = \min_{cd=n}(c+d) = \min_{d\mid n}(d+n/d),
\end{equation*}
where $d\mid n$ means that $d$ divides $n$.

Let $F(n)=n/s(n)$ denote the area-to-semi\-perimeter ratio in which we
are interested. We want to investigate the integers $n$ such that
$F(n)$ is large, and so we define the set $\A$ of ``record-breakers''
for the function $F$ as follows:
\begin{equation}
\A = \{n\in\N\colon F(k)\le F(n) \hbox{ for all }k\le n\}. \label{cdef}
\end{equation}
(Well, the ``record-tiers'' are also included in $\A$.) Then it is
clear after a moment's thought that to solve Grantham's problem for a
given number $N$, we simply need to find the largest element of $\A$
not exceeding $N$.

By computing all possible factorizations of the numbers up to 200 by
brute force, we can make a list of the first 59 elements of $\A$:

\story{$\A=\{$1, 2, 3, 4, 6, 8, 9, 12, 15, 16, 18, 20, 24, 25, 28, 30,
35, 36, 40, 42, 48, 49, 54, 56, 60, 63, 64, 70, 72, 77, 80, 81, 88,
90, 96, 99, 100, 108, 110, 117, 120, 121, 130, 132, 140, 143, 144,
150, 154, 156, 165, 168, 169, 176, 180, 182, 192, 195, 196, \dots\}}

\noindent If we write, in place of the elements $n\in\A$, the
dimensions of the rectangles with area $n$ and least \sp, we
obtain

\story{$\A=\{$1\tms1, 1\tms2, 1\tms3, 2\tms2, 2\tms3, 2\tms4, 3\tms3,
3\tms4, 3\tms5, 4\tms4, 3\tms6, 4\tms5, 4\tms6, 5\tms5, 4\tms7,
5\tms6, 5\tms7, 6\tms6, 5\tms8, 6\tms7, 6\tms8, 7\tms7, 6\tms9,
7\tms8, 6\tms10, 7\tms9, 8\tms8, 7\tms10, 8\tms9, 7\tms11, 8\tms10,
9\tms9, 8\tms11, 9\tms10, 8\tms12, 9\tms11, 10\tms10, 9\tms12,
10\tms11, 9\tms13, 10\tms12, 11\tms11, 10\tms13, 11\tms12, 10\tms14,
11\tms13, 12\tms12, 10\tms15, 11\tms14, 12\tms13, 11\tms15, 12\tms14,
13\tms13, 11\tms16, 12\tms15, 13\tms14, 12\tms16, 13\tms15, 14\tms14,
\dots\},}

\noindent a list that exhibits a tantalizing promise of pattern! The
interested reader is invited to try to determine the precise pattern
of the set $\A$, before reading into the next section in which the
secret will be revealed. One thing we immediately notice, though, is
that the dimensions of each of these rectangles are almost (or
exactly) equal. For this reason, we will call the elements of $\A$
{\it \asqs}.  This supports our intuition about what the answers to
Grantham's problem should be, since after all, Farmer Ted would build
his rectangles with precisely equal sides if he weren't hampered by
the integral policies of (the ironically-named) Builders Square.

From the list of the first 59 \asqs, we find that 182 is the largest
\asq\ not exceeding 190. Therefore, Farmer Ted should build a chicken
coop with area 182 square feet; and indeed, a $13'\tms14'$ coop would
give him about 3.37 square feet of coop space per foot of chicken wire
purchased, which is more cost-effective than the options he thought of
back in Section~\ref{introsec}. But what about next time, when Farmer
Ted wants to build a supercoop on the 8,675,309 square feet of land he
has to spare, or even more? Eventually, computations will need to
give way to a better understanding of $\A$.

Our specific goals in this paper will be to answer the following
questions:

\begin{enumerate}
\item Can we describe $\A$ more explicitly? That is, can we
characterize when a number $n$ is an \asq\ with a description
that refers only to $n$ itself, rather than all the numbers smaller
than $n$? Can we find a formula for the number of \asqs\ not
exceeding a given positive number $x$?
\item Can we quickly compute the largest \asq\ not exceeding $N$, for
a given number $N$? We will describe more specifically what we mean by
``quickly'' in the next section, but for now we simply say that we'll
want to avoid both brute force searches and computations that involve
factoring integers.
\end{enumerate}

\noindent In the next section, we will find that these questions have
surprisingly elegant answers.

\section{Remarkable Results}\label{resultsec}\noindent
Have you uncovered the pattern of the \asqs? One detail you might have
noticed is that all numbers of the form $m\tms m$ and $(m-1)\tms m$,
and also $(m-1)\tms(m+1)$, seem to be \asqs. (If not, maybe we should
come up with a better name for the elements of~$\A$!) This turns out
to be true, as we will see in Lemma \ref{intuitlem} below. The problem
is that there are other \asqs\ than these---3\tms6, 4\tms7, 5\tms8,
6\tms9, 6\tms10---and the ``exceptions'' seem to become more and more
numerous\dots. Even so, it will be convenient to think of the
particular \asqs\ of the form $m\tms m$ and $(m-1)\tms m$ as
``punctuation'' of a sort for~$\A$. To this end, we will define a {\it
flock} to be the set of \asqs\ between $(m-1)^2+1$ and $m(m-1)$, or
between $m(m-1)+1$ and $m^2$, including the endpoints in both cases.

If we group the rectangles corresponding to the \asqs\ into flocks
in this way, indicating the end of each flock by a semicolon, we
obtain:

\story{\newcommand{\bc}{;\ \hskip2pt}$\A=\{$1\tms1\bc 1\tms2\bc
1\tms3, 2\tms2\bc 2\tms3\bc 2\tms4, 3\tms3\bc 3\tms4\bc 3\tms5,
4\tms4\bc 3\tms6, 4\tms5\bc 4\tms6, 5\tms5\bc 4\tms7, 5\tms6\bc
5\tms7, 6\tms6\bc 5\tms8, 6\tms7\bc 6\tms8, 7\tms7\bc 6\tms9,
7\tms8\bc 6\tms10, 7\tms9, 8\tms8\bc 7\tms10, 8\tms9\bc 7\tms11,
8\tms10, 9\tms9\bc 8\tms11, 9\tms10\bc 8\tms12, 9\tms11, 10\tms10\bc
9\tms12, 10\tms11\bc 9\tms13, 10\tms12, 11\tms11\bc 10\tms13,
11\tms12\bc 10\tms14, 11\tms13, 12\tms12\bc 10\tms15, 11\tms14,
12\tms13\bc 11\tms15, 12\tms14, 13\tms13\bc 11\tms16, 12\tms15,
13\tms14\bc 12\tms16, 13\tms15, 14\tms14\bc \dots\}}
\label{seglist}

\noindent It seems that all of the rectangles in a given flock have
the same \sp; this also turns out to be true, as we will see in Lemma
\ref{seqlem} below. The remaining question, then, is to determine
which rectangles of the common \sp\ a given flock contains. At
first it seems that all rectangles of the ``right'' \sp\ will be in
the flock as long as their area exceeds that of the last rectangle
in the preceding flock, but then we note a few omissions---2\tms5,
3\tms7, 4\tms8, 5\tms9, 5\tms10---which also become more
numerous if we extend our computations of~$\A$\dots.

But as it happens, this question can be resolved, and we can actually
determine exactly which numbers are \asqs, as our main theorem
indicates. Recall that $\floor x$ denotes the greatest integer not
exceeding~$x$.

\medskip
\noindent{\bf Main Theorem.} \it
For any integer $m\ge2$, the set of \asqs\ between $(m-1)^2+1$ and
$m^2$ (inclusive) consists of two flocks, the first of which is
\begin{equation*}
%\A\cap( (m-1)^2,m(m-1) ] = 
\{ (m+a_m)(m-a_m-1), (m+a_m-1)(m-a_m), \dots, (m+1)(m-2), m(m-1) \}
\end{equation*}
where $a_m = \floor{(\sqrt{2m-1}-1)/2}$, and the second of which is
\begin{equation*}
%\A\cap( m(m-1),m^2 ] = 
\{ (m+b_m)(m-b_m), (m+b_m-1)(m-b_m+1), \dots, (m+1)(m-1), m^2 \}
\end{equation*}
where $b_m = \floor{\sqrt{m/2}}$.\rm
\medskip

The Main Theorem allows us to easily enumerate the \asqs\ in order,
but if we simply want an explicit characterization of \asqs\ without
regard to their order, there turns out to be one that is extremely
elegant. To describe it, we recall that the {\it triangular numbers\/}
$\{0,1,3,6,10,15,\dots\}$ are the numbers $t_n = {n\choose2} =
n(n-1)/2$ (Conway and Guy \cite{ConGuy} describe many interesting
properties of these and other ``figurate'' numbers). We let $T(x)$
denote the number of triangular numbers not exceeding $x$. (Notice
that in our notation, $t_1 = {1\choose2} = 0$ is the first triangular
number, so that $T(1)=2$, for instance.) Then an alternate
interpretation of the Main Theorem is the following:

\begin{corollary}
The \asqs\ are precisely those integers that can be written in the
form $k(k+h)$, for some integers $k\ge1$ and $0\le h\le T(k)$.
\label{tricor}
\end{corollary}

\begin{figure}[htbf]
\smaller\smaller
\hskip.25in\psfig{figure=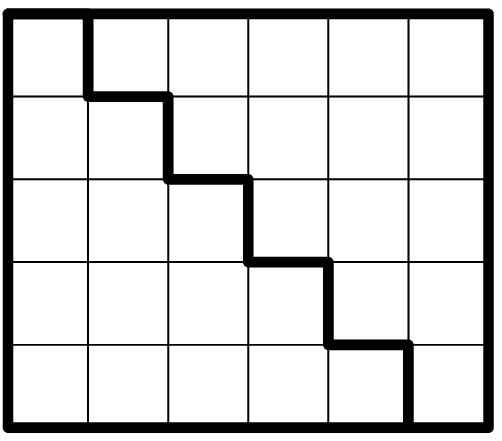}
\put{$m(m-1)=t_m+t_m$}{.08in}{1.4in}
\put{$m^2=t_m+t_{m+1}$}{.8in}{.5in}
\hskip1.9in\psfig{figure=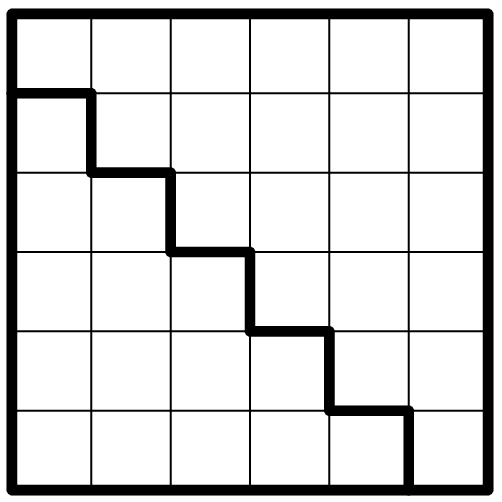}
\caption{Two triangular integers invoke an \asq}
\label{tplustfig}
% 1 grid box = 1/4 inch
\end{figure}

\noindent It is perhaps not so surprising that the triangular numbers
are connected to the \asqs---after all, adding $t_m$ to itself or to
$t_{m+1}$ yields \asqs\ of the form $m(m-1)$ or $m^2$, respectively
(Figure \ref{tplustfig} illustrates this for $m=6$). In any case, the
precision of this characterization is quite attractive and unexpected,
and it is conceivable that Corollary \ref{tricor} has a direct proof
that doesn't use the Main Theorem. We leave this as an open problem
for the reader.

In a different direction, we can use the Main Theorem's precise
enumeration of the \asqs\ in each flock to count the number of \asqs\
quite accurately.

\begin{corollary}
Let $A(x)$ denotes the number of \asqs\ not exceeding $x$. Then for
$x\ge1$,
\begin{equation*}
A(x) = \frac{2\sqrt2}3 x^{3/4} + \frac12x^{1/2} + R(x),
%\label{cxasym}
\end{equation*}
where $R(x)$ is an oscillating term whose order of magnitude is
$x^{1/4}$.
\label{asymcor}
\end{corollary}

\begin{figure}[ht]
\smaller\smaller
\hskip.4cm
\psfig{figure=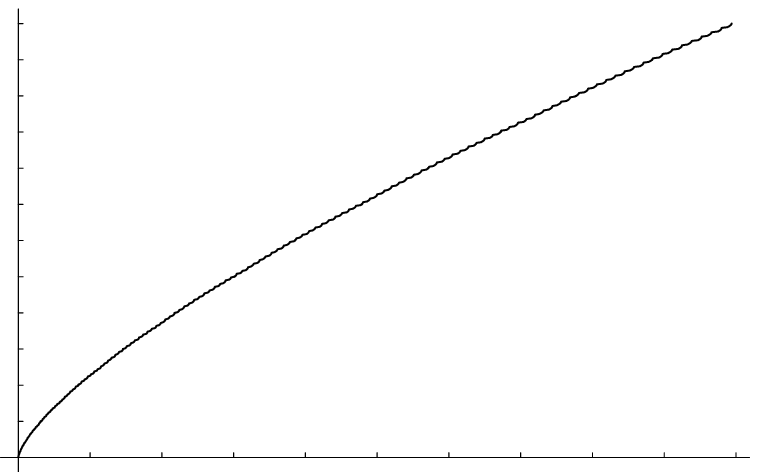}
\put{$x$}{-.2cm}{.4cm}
\put{$A(x)$}{-3.7cm}{3.9cm}
{\smaller\smaller
 \dimen1=-15.25pt\dimen2=-2pt\dimen3=-41.5pt
% \dimen1=-36pt\dimen2=-2pt\dimen3=-41.5pt
 \put{5000}{\dimen1}{\dimen2}
 \advance\dimen1by\dimen3 \put{4000}{\dimen1}{\dimen2}
 \advance\dimen1by\dimen3 \put{3000}{\dimen1}{\dimen2}
 \advance\dimen1by\dimen3 \put{2000}{\dimen1}{\dimen2}
 \advance\dimen1by\dimen3 \put{1000}{\dimen1}{\dimen2}
 \dimen1=-228.7pt\dimen2=23.3pt\dimen3=20.9pt
 \put{100}{\dimen1}{\dimen2}
 \advance\dimen2by\dimen3 \put{200}{\dimen1}{\dimen2}
 \advance\dimen2by\dimen3 \put{300}{\dimen1}{\dimen2}
 \advance\dimen2by\dimen3 \put{400}{\dimen1}{\dimen2}
 \advance\dimen2by\dimen3 \put{500}{\dimen1}{\dimen2}
 \advance\dimen2by\dimen3 \put{600}{\dimen1}{\dimen2}
}
\hskip.3cm
\psfig{figure=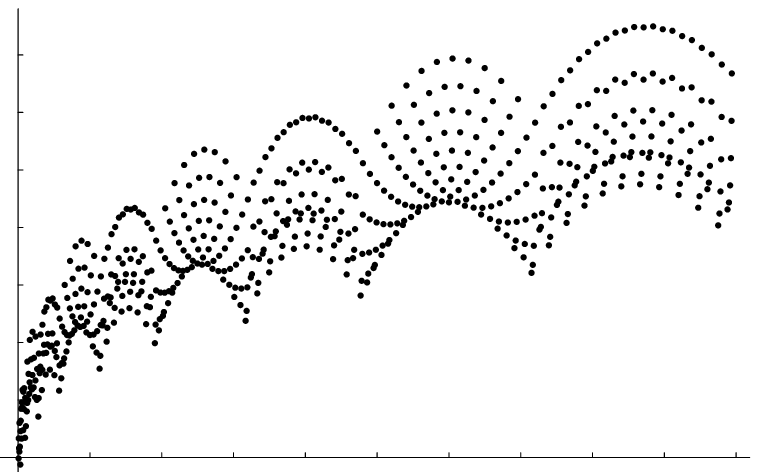}
\put{$x$}{-.2cm}{.4cm}
\put{$R(x)$}{-6cm}{3.9cm}
{\smaller\smaller
 \dimen1=-15.25pt\dimen2=-2pt\dimen3=-41.5pt
 \put{5000}{\dimen1}{\dimen2}
 \advance\dimen1by\dimen3 \put{4000}{\dimen1}{\dimen2}
 \advance\dimen1by\dimen3 \put{3000}{\dimen1}{\dimen2}
 \advance\dimen1by\dimen3 \put{2000}{\dimen1}{\dimen2}
 \advance\dimen1by\dimen3 \put{1000}{\dimen1}{\dimen2}
 \dimen1=-221pt\dimen2=36pt\dimen3=33pt
 \put 2{\dimen1}{\dimen2}
 \advance\dimen2by\dimen3 \put 4{\dimen1}{\dimen2}
 \advance\dimen2by\dimen3 \put 6{\dimen1}{\dimen2}
}
\caption{Superficial steadiness of $A(x)$, mesmerizing meanderings of
$R(x)$}
\label{axrxfig}
\end{figure}

\noindent A graph of $A(x)$ (see Figure \ref{axrxfig}) exhibits a
steady growth with a little bit of a wiggle. When we isolate $R(x)$ by
subtracting the main term $2\sqrt2x^{3/4}/3+x^{1/2}/2$ from $A(x)$,
the resulting graph (Figure \ref{axrxfig}, where we have plotted a
point every time $x$ passes an \asq) is a pyrotechnic, almost
whimsical display that seems to suggest that our computer code needs
to be rechecked. Yet this is the true nature of $R(x)$. When we prove
Corollary \ref{asymcor} (in a more specific and precise form) in
Section \ref{cacsec}, we will see that there are two reasons that the
``remainder term'' $R(x)$ oscillates: there are oscillations on a
local scale because the \asqs\ flock towards the right half of each
interval of the form $((m-1)^2,m(m-1)]$ or $(m(m-1),m^2]$, and
oscillations on a larger scale for a less obvious reason.

These theoretical results about the structure of the \asqs\ address
question 1 nicely, and we turn our attention to the focus of question
2, the practicality of actually computing answers to questions about
\asqs. Even simple tasks like printing out a number or adding two
numbers together obviously take time for a computer to perform, and
they take longer for bigger numbers. To measure how the computing time
needed for a particular computation increases as the size of the input
grows, let $f(k)$ denote the amount of time it takes to perform the
calculation on a $k$-digit number. Of course, the time could depend
significantly on which $k$-digit number we choose; what we mean is the
worst-case scenario, so that the processing time is at most $f(k)$ no
matter which $k$-digit number we choose.

We say that a computation runs in {\it polynomial time\/} if this
function $f(k)$ grows only as fast as a polynomial in $k$, i.e., if
there are positive constants $A$ and $B$ such that
$f(k)<Ak^B$. Generally speaking, the computations that we consider
efficient to perform on very large inputs are those that run in
polynomial time. (Because we are only concerned with this category of
computations as a whole, it doesn't matter if we write our numbers in
base 10 or base 2 or any other base, since this only multiplies the
number of digits by a constant factor like $\log_210$.)

All of our familiar arithmetic operations $+$, $-$, $\tms$, $\div$,
$\sqrt\cdot$, $\floor\cdot$ and so on have polynomial-time
algorithms. On the other hand, performing a calculation on each of the
numbers from 1 to the input $n$, or even from 1 to $\sqrt n$, etc.,~is
definitely not polynomial-time. Thus computing \asqs\ by their
definition, which involves comparing $F(n)$ with all of the preceding
$F(k)$, is not efficient for large $n$. Furthermore, the obvious
method of factoring numbers---testing all possible divisors in
turn---is not polynomial-time for the same reason. While there are
faster ways to factor numbers, at this time there is still no known
polynomial-time algorithm for factoring numbers; so even factoring a
single number would make an algorithm inefficient. (Dewdney
\cite{dewdney} writes about many facets of algorithms, including this
property of running in polynomial time, while Pomerance \cite{pom}
gives a more detailed discussion of factoring algorithms and their
computational complexity.)

Fortunately, the Main Theorem provides a way to compute \asqs\ that
avoids both factorization and brute-force enumeration. In fact, we can
show that all sorts of computations involving \asqs\ are efficient:

\begin{corollary}
There are polynomial-time algorithms to perform each of the following
tasks, given a positive integer $N$:
\begin{enumerate}
\item[{\rm(a)}]determine whether $N$ is an \asq, and if so determine
the dimensions of the optimal rectangle;
\item[{\rm(b)}]find the greatest \asq\ not exceeding $N$, including
the dimensions of the optimal rectangle;
\item[{\rm(c)}]compute the number $A(N)$ of \asqs\ not exceeding $N$;
\item[{\rm(d)}]find the $N$th \asq, including the dimensions of the
optimal rectangle.
\end{enumerate}
\label{polycor}
\end{corollary}

\noindent We reiterate that these algorithms work without ever
factoring a single integer. Corollary~\ref{polycor}, together with our
lack of a polynomial-time factoring algorithm, has a rather
interesting implication: for large values of $N$, it is much faster to
compute the most cost-effective chicken coop (in terms of
area-to-semiperimeter ratio) with area {\it at most\/} $N$ than it is
to compute the most cost-effective chicken coop with area {\it equal
to\/} $N$, a somewhat paradoxical state of affairs! Nobody ever said
farming was easy\dots.

\section{The Theorem Thought Through}\label{ththsec}\noindent
Before proving the Main Theorem, we need to build up a stockpile
of easy lemmas. The first of these simply confirms our expectations
that the most cost-effective rectangle of a given area is the one
whose side lengths are as close together as possible, and also
provides some inequalities for the functions $s(n)$ and $F(n)$. Let us
define $d(n)$ to be the largest divisor of $n$ not exceeding $\sqrt n$
and $d'(n)$ the smallest divisor of $n$ that is at least $\sqrt n$, so
that $d'(n)=n/d(n)$.

\begin{lemma}
The rectangle with integer side lengths and area $n$ that has the
smallest \sp\ is the one with dimensions $d(n)\tms
d'(n)$. In other words,
\begin{equation*}
s(n) = d(n) + d'(n).
\end{equation*}
We also have the inequalities
\begin{equation*}
s(n) \ge 2\sqrt n \quad\hbox{and}\quad F(n) \le \sqrt n/2.
\end{equation*}
\label{firstlem}
\end{lemma}

\begin{proof}
For a fixed positive number $n$, the function $f(t)=t+n/t$ has
derivative $f'(t)=1-n/t^2$ which is negative for $1\le t<\sqrt n$, and
therefore $t+n/t$ is a decreasing function of $t$ in that range. Thus
if we restrict our attention to those $t$ such that both $t$ and $n/t$
are positive integers (in other words, $t$ is an integer dividing
$n$), we see that the expression $t+n/t$ is minimized when
$t=d(n)$. We therefore have
\begin{equation*}
s(n)=d(n)+n/d(n)\ge2\sqrt n,
\end{equation*}
where the last inequality follows from the Arithmetic Mean/Geometric
Mean inequality; and the inequality for $F(n)$ then follows directly
from the definition of $F$.
\end{proof}

If we have a number $n$ written as $c\tms d$ where $c$ and $d$ are
pretty close to $\sqrt n$, when can we say that there isn't some
better factorization out there, so that $s(n)$ is really equal to
$c+d$? The following lemma gives us a useful criterion.

\begin{lemma}
If a number $n$ satisfying $(m-1)^2<n\le m(m-1)$ has the form
$n=(m-a-1)(m+a)$ for some number $a$, then $s(n)=2m-1$, and
$d(n)=m-a-1$ and $d'(n)=m+a$. Similarly, if a number $n$ satisfying
$m(m-1)<n\le m^2$ has the form $n=m^2-b^2$ for some number $b$, then
$s(n)=2m$, and $d(n)=m-b$ and $d'(n)=m+b$.
\label{nohidelem}
\end{lemma}

\begin{proof}
First let's recall that for any positive real numbers $\alpha$ and
$\beta$, the pair of equations $r+s=\alpha$ and $rs=\beta$ have a
unique solution $(r,s)$ with $r\le s$, as long as the Arithmetic
Mean/Geometric Mean inequality $\alpha/2\ge\sqrt\beta$ holds. This is
because $r$ and $s$ will be the roots of the quadratic polynomial
$t^2-\alpha t+\beta$, which has real roots when its discriminant
$\alpha^2-4\beta$ is nonnegative, i.e., when $\alpha/2\ge\sqrt\beta$.

Now if $n=(m-a-1)(m+a)$, then clearly $s(n)\le(m-1-a)+(m+a)=2m-1$ by
the definition of $s(n)$. On the other hand, by Lemma \ref{firstlem}
we know that $s(n)\ge2\sqrt n>2(m-1)$, and so $s(n)=2m-1$
exactly. We now know that
$d(n)d'(n)=n=(m-a-1)(m+a)$ and
\begin{equation*}
d(n)+d'(n) = s(n) = 2m-1 = (m-a-1)+(m+a),
\end{equation*}
and of course $d(n)\le d'(n)$ as well; by the argument of the previous
paragraph, we conclude that $d(n)=m-a-1$ and $d'(n)=m+a$. This
establishes the first assertion of the lemma, and a similar argument
holds for the second assertion.
\end{proof}

Of course, if a number $n$ satisfies $s(n)=2m$ for some $m$, then $n$
can be written as $n=cd$ with $c\le d$ and $c+d=2m$; and letting
$b=d-m$, we see that $n=cd=(2m-d)d=(m-b)(m+b)$. A similar statement is
true if $s(n)=2m-1$, and so we see that the converse of Lemma
\ref{nohidelem} also holds. We also remark that in the statement of
the lemma, the two expressions $m(m-1)$ can be replaced by
$(m-1/2)^2=m(m-1)+1/4$ if we wish.

Lemma \ref{nohidelem} implies in particular that for $m\ge2$,
\begin{equation*}
s(m^2)=2m, \quad s(m(m-1))=2m-1, \quad\hbox{and } s((m-1)(m+1))=2m,
\end{equation*}
and so
\begin{equation*}
F(m^2) = \frac m2, \quad F(m^2-m) = {m(m-1)\over2m-1}, \quad\hbox{and
} F(m^2-1) = \frac{m^2-1}{2m}.
\end{equation*}
Using these facts, we can verify our theory that these numbers are
always \asqs.

\begin{lemma}
Each positive integer of the form $m^2$, $m(m-1)$, or $m^2-1$ is an
\asq.\label{intuitlem}
\end{lemma}

It is interesting to note that these are precisely those integers $n$
that are divisible by $\floor{\sqrt n}$ (see \cite{sqrtn}), one of
the many interesting things that can be discovered by referring to
Sloane and Plouffe~\cite{sloplo}.
\medskip

\begin{proof}
We verify directly that such numbers satisfy the condition in the
definition (\ref{cdef}) of~$\A$. If $k<m^2$, then by Lemma
\ref{firstlem} we have $F(k)\le\sqrt{k}/2 < m/2=F(m^2)$, and so $m^2$
is an \asq. Similarly, if $k<m(m-1)$, then again
\begin{equation*}
F(k) \le \frac{\sqrt{k}}2 \le \frac{\sqrt{m^2-m-1}}2 <
{m(m-1)\over2m-1} = F(m(m-1)),
\end{equation*}
where the strict inequality can be verified as a ``fun'' algebraic
exercise. Thus $m(m-1)$ is also an \asq. A similar argument shows that
$m^2-1$ is also an \asq.
\end{proof}

Now we're getting somewhere! Next we show that the \sps\ of the
rectangles corresponding to the \asqs\ in a given flock are all
equal, as we observed at the beginning of Section \ref{resultsec}.

\begin{lemma}
Let $m\ge2$ be an integer. If $n$ is an \asq\ satisfying
$(m-1)^2<n\le m(m-1)$, then $s(n)=2m-1$; similarly, if $n$ is an \asq\
satisfying $m(m-1)<n\le m^2$, then $s(n)=2m$.
\label{seqlem}
\end{lemma}

\begin{proof}
If $n=m(m-1)$, we have already shown that $s(n)=2m-1$. If $n$
satisfies $(m-1)^2<n<m(m-1)$, then by Lemma \ref{firstlem} we have
$s(n)\ge2\sqrt n>2(m-1)$. On the other hand, since $n$ is an \asq\
exceeding $(m-1)^2$, we have
\begin{equation*}
\frac{m-1}2 = F((m-1)^2) \le F(n) = \frac n{s(n)} < {m(m-1)\over s(n)},
\end{equation*}
and so $s(n)<2m$. Therefore $s(n)=2m-1$ in this case.

Similarly, if $n$ satisfies $m(m-1)<n<m^2$, then $s(n)\ge2\sqrt
n\ge2\sqrt{m^2-m+1} > 2m-1$; on the other hand,
\begin{equation*}
{m(m-1)\over2m-1} = F(m(m-1)) \le F(n) = \frac n{s(n)} \le {m^2-1\over
s(n)},
\end{equation*}
and so $s(n)<(m+1)(2m-1)/m<2m+1$. Therefore $s(n)=2m$ in this case.
\end{proof}

Finally, we need to exhibit some properties of the sequences
$a_m$ and $b_m$ defined in the statement of the Main Theorem.

\begin{lemma}
Define $a_m=\floor{(\sqrt{2m-1}-1)/2}$ and
$b_m=\floor{\sqrt{m/2}}$. For any integer $m\ge2$:
\begin{enumerate}
\item[{\rm(a)}]$a_m\le b_m\le a_m+1$;
\item[{\rm(b)}]$b_m=\floor{m/\sqrt{2m-1}}$; and
\item[{\rm(c)}]$a_m+b_m=\floor{\sqrt{2m}}-1$.
\end{enumerate}
\label{floorlem}
\end{lemma}

We omit the proof of this lemma since it is tedious but
straightforward. The idea is to show that in the sequences $a_m$,
$b_m$, $\floor{m/\sqrt{2m-1}}$, and so on, two consecutive terms are
either equal or else differ by 1, and then to determine precisely for
what values of $m$ the differences of 1 occur.

Armed with these lemmas, we are now ready to furnish a proof of
the Main Theorem.\medskip

\begin{pflike}{Proof of the Main Theorem:}
Fix an integer $m\ge2$. By Lemma \ref{seqlem}, every \asq\ $n$ with
$(m-1)^2<n\le m(m-1)$ satisfies $s(n)=2m-1$; while by Lemma
\ref{nohidelem}, the integers $(m-1)^2<n\le m(m-1)$ satisfying
$s(n)=2m-1$ are precisely the elements of the form $n_a=(m-a-1)(m+a)$
that lie in that interval. Thus it suffices to determine which of the
$n_a$ are \asqs.

Furthermore, suppose that $n_a$ is an \asq\ for some $a\ge1$. Then
$F(n_a)\ge F(n)$ for all $n<n_a$ by the definition of $\A$, while
$F(n_a)>F(n)$ for all $n_a<n<n_{a-1}$ since we've already concluded
that no such $n$ can be an \asq. Moreover, $n_{a-1}>n_a$ and
$s(n_{a-1})=2m-1=s(n_a)$, so $F(n_{a-1})>F(n_a)$, and thus $n_{a-1}$
is an \asq\ as well. Therefore it suffices to find the largest value
of $a$ (corresponding to the smallest $n_a$) such that $n_a$ is an
\asq.

By Lemma \ref{intuitlem}, we know that $(m-1)^2$ is an \asq, and so we
need to find the largest $a$ such that $F(n_a)\ge F((m-1)^2)$, i.e.,
\begin{equation*}
{(m-a-1)(m+a)\over2m-1} \ge \frac{m-1}2,
\end{equation*}
which is the same as $2a(a+1)+1\le m$. By completing the square and
solving for $a$, we find that this inequality is equivalent to
\begin{equation}
\frac{-\sqrt{2m-1}-1}2 \le a \le \frac{\sqrt{2m-1}-1}2,  \label{aineq}
\end{equation}
and so the largest integer $a$ satisfying the inequality is exactly
$a=\floor{(\sqrt{2m-1}-1)/2}=a_m$, as defined in the statement of
the Main Theorem. This establishes the first part of the theorem.

By the same reasoning, it suffices to find the largest value of $b$
such that $F(m^2-b^2)\ge F(m(m-1))$, i.e.,
\begin{equation*}
{m^2-b^2\over2m} \ge {m(m-1)\over2m-1},
\end{equation*}
which is the same as
\begin{equation}
b^2\le m^2/(2m-1)  \label{bineq}
\end{equation}
or $b\le\floor{m/\sqrt{2m-1}}$. But by Lemma \ref{floorlem}(b),
$\floor{m/\sqrt{2m-1}}=b_m$ for $m\ge2$, and so the second part of the
theorem is established.
\end{pflike}

With the Main Theorem now proven, we remark that Lemma
\ref{floorlem}(c) implies that for any integer $m\ge2$, the number of
\asqs\ in the two flocks between $(m-1)^2+1$ and $m^2$ is exactly
$(1+a_m)+(1+b_m) = 1+\floor{\sqrt{2m}}$, while Lemma \ref{floorlem}(a)
implies that there are either equally many in the two flocks or else
one more in the second flock than in the first.

\section{Taking Notice of Triangular Numbers}\label{tntnsec}\noindent
Our next goal is to derive Corollary \ref{tricor} from the Main
Theorem. First we establish a quick lemma giving a closed-form
expression for $T(x)$, the number of triangular numbers not
exceeding~$x$.

\begin{lemma}
For all $x\ge0$, we have $T(x)=\floor{\sqrt{2x+1/4}+1/2}$.  \label{Txlem}
\end{lemma}

\begin{proof}
$T(x)$ is the number of positive integers $n$ such that $t_n\le x$, or
$n(n-1)/2 \le x$. This inequality is equivalent to $(n-1/2)^2 \le
2x+1/4$, or $-\sqrt{2x+1/4}+1/2 \le n \le \sqrt{2x+1/4}+1/2$. The
left-hand expression never exceeds $1/2$, and so $T(x)$ is simply the
number of positive integers $n$ such that $n \le \sqrt{2x+1/4}+1/2$;
in other words, $T(x) = \floor{\sqrt{2x+1/4}+1/2}$ as desired.
\end{proof}

\medskip
\begin{pflike}{Proof of Corollary \ref{tricor}:}
Suppose first that $n=k(k+h)$ for some integers $k\ge1$ and $h\le
T(k)$. Let $k'=k+h$, and define
\begin{equation*}
\oecases{m=k+(h+1)/2 \and a=(h-1)/2}{m=k+h/2 \and b=h/2},
\end{equation*}
so that
\begin{equation*}
\oecases{k=m-a-1 \and k'=m+a}{k=m-b \and k'=m+b}.
\end{equation*}
We claim that
\begin{equation}
\oecases{(m-1)^2 < (m-a-1)(m+a) \le (m-1/2)^2}{(m-1/2)^2 < m^2-b^2 \le
m^2}.  \label{segclaim}
\end{equation}
To see this, note that in terms of $k$ and $h$, these inequalities
become
\begin{equation*}
\big( k+\frac{h-1}2 \big)^2 < k(k+h) \le \big( k+\frac h2 \big)^2.
%\label{khclaim}
\end{equation*}
A little bit of algebra reveals that the right-hand inequality is
trivially satisfied while the left-hand inequality is true provided
that $h<2\sqrt k+1$. However, from Lemma \ref{Txlem} we see that
\begin{equation*}
T(k) = \floor{\sqrt{2k+1/4}+1/2} \le \sqrt{2k+1/4}+1/2 < 2\sqrt k+1
\end{equation*}
for $k\ge1$. Since we are assuming that $h\le T(k)$, this shows that
the inequalities (\ref{segclaim}) do indeed hold.

Because of these inequalities, we may apply Lemma \ref{nohidelem} (see
the remarks following the proof of the lemma) and conclude that
\begin{equation*}
\oecases{s(n)=2m-1,\, d(n)=m-a-1, \and d'(n)=m+a}{s(n)=2m,\, d(n)=m-b,
\and d'(n)=m+b}.
\end{equation*}
Consequently, the Main Theorem asserts that $n$ is an \asq\ if and
only if
\begin{equation}
\oecases{a\le a_m}{b\le b_m},  \label{abbd}
\end{equation}which by the definitions of $a$, $b$, and $m$ is the same as
\begin{equation*}
\oecases{(h-1)/2 \le \floor{(\sqrt{2m-1}-1)/2} =
\floor{(\sqrt{2k+h}-1)/2}}{h/2 \le \floor{\sqrt{m/2}} =
\floor{\sqrt{k/2+h/4}}}.
\end{equation*}
Since in either case, the left-hand side is an integer, the
greatest-integer brackets can be removed from the right-hand side,
whence both cases reduce to $h\le\sqrt{2k+h}$. From here, more algebra
reveals that this inequality is equivalent to $h\le\sqrt{2k+1/4}+1/2$;
and since $h$ is an integer, we can add greatest-integer brackets to
the right-hand side, thus showing that the inequality (\ref{abbd}) is
equivalent to $h\le T(k)$ (again using Lemma \ref{Txlem}). In
particular, $n$ is indeed an \asq.

This establishes one half of the characterization asserted by
Corollary \ref{tricor}. Conversely, suppose we are given an \asq\ $n$,
which we can suppose to be greater than 1 since 1 can obviously be
written as $1(1+0)$. If we let $h=d'(n)-d(n)$, then the Main Theorem
tells us that
\begin{equation*}
\oecases{n=(m-a-1)(m+a),\, d(n)=m-a-1, \and d'(n)=m+a}{n=m^2-b^2,\,
d(n)=m-b, \and d'(n)=m+b}{}
\end{equation*}
for some integers $m\ge2$ and either $a$ with $0\le a\le a_m$ or $b$
with $0\le b\le b_m$. If we set $k=d(n)$, then certainly
$n=k(k+h)$. Moreover, the algebraic steps showing that the inequality
(\ref{abbd}) is equivalent to $h\le T(k)$ are all reversible; and
(\ref{abbd}) does in fact hold, since we are assuming that $n$ is an
\asq. Therefore $n$ does indeed have a representation of the form
$k(k+h)$ with $0\le h\le T(k)$. This establishes the corollary.
\end{pflike}

We take a slight detour at this point to single out some special
\asqs. Let us make the convention that the $k$th flock refers to the
flock of \asqs\ with \sp\ $k$, so that the first flock is actually
empty, the second and third poor flocks contain only $1=1\tms1$ and
$2=1\tms2$, respectively, the fourth flock contains $3=1\tms3$ and
$4=2\tms2$, and so on. The Main Theorem tells us that $a_m$ and $b_m$
control the number of \asqs\ in the odd-numbered and even-numbered
flocks, respectively; thus every so often, a flock will have one more
\asq\ than the preceding flock of the same ``parity''. We'll let a
{\it pioneer\/} be an \asq\ that begins one of these suddenly-longer
flocks.

For instance, from the division of $\A$ into flocks on page
\ref{seglist}, we see that the 4th flock \{1\tms3, 2\tms2\} is
longer than the preceding even-numbered flock \{1\tms1\}, so
$1\tms3=3$ is the first pioneer; the 9th flock \{3\tms6, 4\tms5\} is
longer than the preceding odd-numbered flock \{3\tms4\}, so
$3\tms6=18$ is the second pioneer; and so on, the next two pioneers
being $6\tms10$ in the 16th flock and $10\tms15$ in the 25th
flock. Now if this isn't a pattern waiting for a proof, nothing is!
The following lemma shows another elegant connection between the
\asqs\ and the squares and triangular numbers.

\begin{corollary}
For any positive integer $j$, the $j$th pioneer equals $t_{j+1}\tms
t_{j+2}$ (where $t_i$ is the $i$th triangular number), which begins
the $(j+1)^2$-th flock. Furthermore, the ``record-tying'' \asqs\
(those whose $F$-values are equal to the $F$-values of their immediate
predecessors in $\A$) are precisely the even-numbered pioneers.
\label{piocor}
\end{corollary}

\begin{proof}
First, Lemma \ref{floorlem}(a) tells us that the odd- and
even-numbered flocks undergo their length increases in alternation, so
that the pioneers alternately appear in the flocks of each parity. The
first pioneer $3=1\tms3$ appears in the 4th flock, and corresponds to
$m=2$ and the first appearance of $b_m=1$ in the notation of the Main
Theorem. Thus the $(2k-1)$-st pioneer will equal $m^2-k^2$, where $m$
corresponds to the first appearance of $b_m=k$. It is easy to see that
the first appearance of $b_m=k$ occurs when $m=2k^2$, in which case
the $(2k-1)$-st pioneer is
\begin{equation*}
m^2-k^2 = (2k^2)^2-k^2 = (2k^2-k)(2k^2+k) = \frac{2k(2k-1)}2
\frac{(2k+1)2k}2 = t_{2k}t_{2k+1}.
\end{equation*}
Moreover, the flock in which this pioneer appears is the $2m$-th or
$(2k)^2$-th flock.

Similarly, the $2k$-th pioneer will equal $(m-k-1)(m+k)$, where $m$
corresponds to the first appearance of $a_m=k$. Again one can show
that the first appearance of $a_m=k$ occurs when $m=2k^2+2k+1$, in
which case the $2k$-th pioneer is
\begin{equation*}
(m-k-1)(m+k) = (2k^2+k)(2k^2+3k+1) = \frac{(2k+1)2k}2
\frac{(2k+2)(2k+1)}2 = t_{2k+1}t_{2k+2}.
\end{equation*}
Moreover, the flock in which this pioneer appears is the $(2m-1)$-st
or $(2k+1)^2$-th flock. This establishes the first assertion of the
corollary.

Since the $F$-values of the \asqs\ form a nondecreasing sequence by
the definition of \asq, to look for \asqs\ with equal $F$-values we
only need to examine consecutive \asqs. Furthermore, two consecutive
\asqs\ in the same flock never have equal $F$-values, since they are
distinct numbers but by Lemma \ref{seqlem} their \sps\ are the
same. Therefore we only need to determine when the last \asq\ in a
flock can have the same $F$-value as the first \asq\ in the
following flock.

The relationship between the $F$-values of these pairs of \asqs\ was
determined in the proof of the Main Theorem. Specifically, the
equality $F((m-1)^2)=F((m-a-1)(m+a))$ holds if and only if the
right-hand inequality in (\ref{aineq}) is actually an equality; this
happens precisely when $m=2a^2+2a+1$, which corresponds to the
even-numbered pioneers as was determined above. On the other hand, the
equality $F(m(m-1))=F(m^2-b^2)$ holds if and only if the inequality
(\ref{bineq}) is actually an equality; but $m^2$ and $2m-1$ are always
relatively prime (any prime factor of $m^2$ must divide $m$ and thus
divides into $2m-1$ with a ``remainder'' of $-1$), implying that
$m^2/(2m-1)$ is never an integer for $m\ge2$, and so the inequality
(\ref{bineq}) can never be an equality. This establishes the second
assertion of the corollary.
\end{proof}

\begin{figure}[htbf]
\smaller\smaller
\centerline{
\psfig{figure=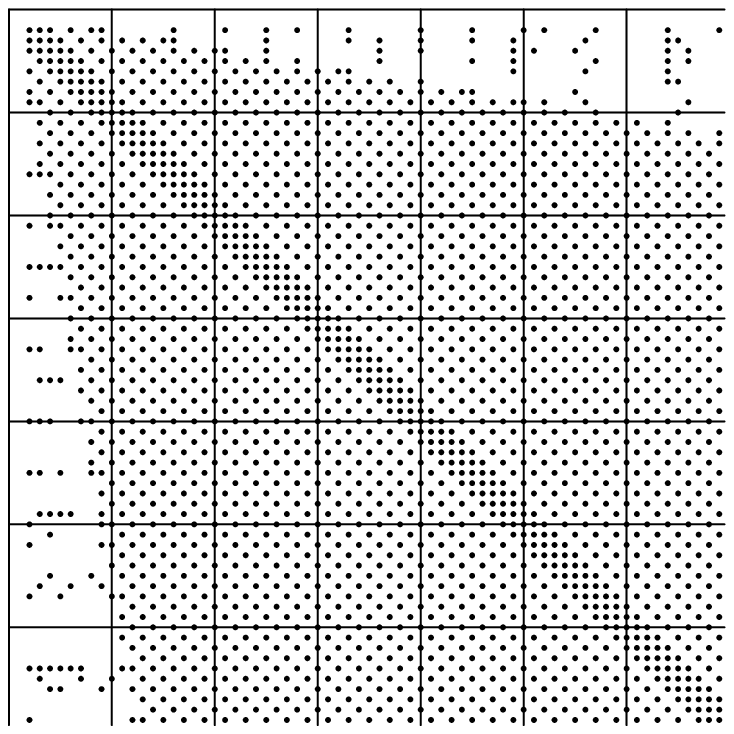}
\put{$m$}{-1.58in}{3.05in}
\put{$n$}{-3.2in}{1.48in}
{
 \smaller\smaller
 \dimen1=-2.255in\dimen2=2.96in\dimen3=.825in
 \put{20}{\dimen1}{\dimen2}
 \advance\dimen1by\dimen3 \put{40}{\dimen1}{\dimen2}
 \advance\dimen1by\dimen3 \put{60}{\dimen1}{\dimen2}
 \dimen1=-3.15in\dimen2=2.083in\dimen3=-.825in
 \put{20}{\dimen1}{\dimen2}
 \advance\dimen2by\dimen3 \put{40}{\dimen1}{\dimen2}
 \advance\dimen2by\dimen3 \put{60}{\dimen1}{\dimen2}
}}
\caption{Amazing \asq\ patterns in products of two triangles}
\label{tmxtnfig}
\end{figure}

We know that all squares are \asqs, and so $t_j^2$ is certainly an
\asq\ for any triangular number $t_j$; also, Corollary \ref{piocor}
tells us that the product $t_jt_{j+1}$ of two consecutive triangular
numbers is always an \asq. This led the author to wonder which numbers
of the form $t_mt_n$ are \asqs. If $m$ and $n$ differ by more than 1,
it would seem that the rectangle of dimensions $t_m\tms t_n$ is not
the most cost-effective rectangle of area $t_mt_n$, and so the author
expected that these products of two triangular numbers would behave
randomly with respect to being \asqs---that is, a few of them might be
but most of them wouldn't. After some computations, however, Figure
\ref{tmxtnfig} emerged, where a point has been plotted in the $(m,n)$
position if and only if $t_mt_n$ is an \asq; and the table exhibited a
totally unexpected regularity.

Of course the symmetry of the table across the main diagonal is to be
expected since $t_mt_n=t_nt_m$. The main diagonal and the first
off-diagonals are filled with plotted points, corresponding to the
\asqs\ $t_m^2$ and $t_mt_{m+1}$; and in hindsight, the second
off-diagonals correspond to
\begin{equation*}
t_mt_{m+2} = \frac{m(m-1)}2\frac{(m+2)(m+1)}2 =
\frac{m^2+m-2}2 \, \frac{m^2+m}2,
\end{equation*}
which is the product of two consecutive integers (since $m^2+m$ is
always even) and is thus an \asq\ as well. But apart from these
central diagonals and some garbage along the edges of the table where
$m$ and $n$ are quite different in size, the checkerboard-like pattern
in the kite-shaped region of the table seems to be telling us
that the only thing that matters in determining whether $t_mt_n$ is an
\asq\ is whether $m$ and $n$ have the same parity!

Once this phenomenon had been discovered, it turned out that the
following corollary could be derived from the prior results in this
paper. We leave the proof of this corollary as a challenge to the
reader.

\begin{corollary}
Let $m$ and $n$ be positive integers with $n-1>m>3n-\sqrt{8n(n-1)}-1$.
Then $t_mt_n$ is an \asq\ if and only if $n-m$ is even.
\end{corollary}

We remark that the function $3n-\sqrt{8n(n-1)}-1$ is asymptotic to
$(3-2\sqrt2)n+(\sqrt2-1)$, which explains the straight lines of slope
$-(3-2\sqrt2)\approx{-0.17}$ and $-1/(3-2\sqrt2)\approx{-5.83}$ that
seem to the eye to separate the orderly central region in Figure
\ref{tmxtnfig} from the garbage along the edges.

\section{Counting and Computing}\label{cacsec}\noindent
In this section we establish Corollaries \ref{asymcor} and
\ref{polycor}. We begin by defining a function $B(x)$ that will serve
as the backbone of our investigation of the \asq\ counting function
$A(x)$. Let $\{x\}=x-\floor x$ denote the fractional part of $x$, and
define the quantities $\gamma=\gamma(x) = \{\sqrt2x^{1/4}\}$ and
$\delta=\delta(x) = \{x^{1/4}/\sqrt2\}$. Let $B(x)=B_0(x)+B_1(x)$,
where
\begin{equation}
B_0(x) = \frac{2\sqrt2}3x^{3/4} + \frac12x^{1/2} + \big(
\frac{2\sqrt2}3 + \frac{\gamma(1-\gamma)}{\sqrt2} \big) x^{1/4}
\label{b0def}
\end{equation}
and
\begin{equation*}
B_1(x) = \frac{\gamma^3}6 - \frac{\gamma^2}4 - \frac{5\gamma}{12}
- \frac\delta2 - 1.
%\label{b1def}
\end{equation*}
We remark that $\gamma=\{2\delta\}$ and that $B_1(x^4)$ is a periodic
function of $x$ with period $\sqrt2$, and so it is easy to check that
the inequalities $-2\le B_1(x)\le-1$ always hold. The following lemma
shows how the strange function $B(x)$ arises in connection with the
\asqs.
\begin{lemma}
For any integer $M\ge1$, we have $A(M^2)=B(M^2)$.  \label{am2bm2}
\end{lemma}

\begin{proof}
As remarked at the end of Section~\ref{ththsec}, the number of
\asqs\ between $(m-1)^2+1$ and $m^2$ is $\floor{\sqrt{2m}}+1$ for
$m\ge2$. Therefore
\begin{equation*}
A(M^2) = 1 + \sum_{m=2}^M (\floor{\sqrt{2m}}+1) = M-1 + \sum_{m=1}^M
\floor{\sqrt{2m}}.
\end{equation*}
It's almost always a good idea to interchange orders of summation
whenever possible---and if there aren't enough summation signs, find a
way to create some more! In this case, we convert the greatest-integer
function into a sum of 1s over the appropriate range of integers:
\begin{equation*}
\begin{split}
A(M^2) &= M-1 + \sum_{m=1}^M \sum_{1\le k\le\sqrt{2m}} 1 \\
&= M-1 + \sum_{1\le k\le\sqrt{2M}} \sum_{k^2/2\le m\le M} 1 \\
&= M-1 + \sum \begin{Sb}1\le k\le\sqrt{2M} \\ k\text{ odd}\end{Sb}
\big( M-\frac{k^2-1}2 \big) + \sum \begin{Sb}1\le k\le\sqrt{2M} \\
k\text{ even}\end{Sb} \big( M-\big( \frac{k^2}2-1 \big) \big).
\end{split}
\end{equation*}
If we temporarily write $\mu$ for $\floor{\sqrt{2M}}$, then
\begin{equation}
\begin{split}
A(M^2) &= M-1 + \mu\big( M+\frac12 \big) - \frac12 \sum_{k=1}^\mu k^2
+ \sum \begin{Sb}k=1 \\ k\text{ even}\end{Sb}^\mu \frac12 \\
&= M(\mu+1) + \frac\mu2-1 - \frac12\frac{\mu(\mu+1)(2\mu+1)}6 -
\frac12\bfloor{\frac\mu2},
\end{split}
\label{muform}
\end{equation}
using the well-known formula for the sum of the first $\mu$ squares.
Since $\floor{\floor x/n} = \floor{x/n}$ for any real number $x\ge0$
and any positive integer $n$, the last term can be written as
\begin{equation*}
\frac12\bfloor{\frac\mu2} = \frac12\bbfloor{\frac{\floor{\sqrt{2M}}}2}
= \frac12\bbfloor{\sqrt{\frac M2}} = \frac12\sqrt{\frac M2} -
\frac12\bigg\{\sqrt{\frac M2}\bigg\} = \frac12\sqrt{\frac M2} -
\frac{\delta(M^2)}2,
\end{equation*}
while we can replace the other occurrences of $\mu$ in equation
(\ref{muform}) by $\sqrt{2M}-\{\sqrt{2M}\} =
\sqrt{2M}-\gamma(M^2)$. Writing $\gamma$ for $\gamma(M^2)$ and
$\delta$ for $\delta(M^2)$, we see that
\begin{equation*}
\begin{split}
A(M^2) &= M(\sqrt{2M}-\gamma+1) + \sqrt{2M}-\gamma-1 \\
&\qquad {}- \frac12
\frac{(\sqrt{2M}-\gamma)(\sqrt{2M}-\gamma+1)(2(\sqrt{2M}-\gamma)+1)}6
- \frac12\sqrt{\frac M2} + \frac\delta2 \\
&= \frac{2\sqrt2}3M^{3/2} + \frac M2 + \big( \frac{2\sqrt2}3 -
\frac{\gamma(1-\gamma)}{\sqrt2} \big) \sqrt M + \frac{\gamma^3}6 -
\frac{\gamma^2}4 - \frac{5\gamma}{12} - 1 - \frac\delta2 = B(M^2)
\end{split}
%\label{bm2eq}
\end{equation*}
after much algebraic simplification. This establishes the lemma.
\end{proof}

Now $B(x)$ is a rather complicated function of $x$, but the next lemma
gives us a couple of ways to predict the behavior of $B(x)$. First, it
tells us how to predict $B(x+y)$ from $B(x)$ if $y$ is small compared
to $x$ (roughly speaking, their difference will be $y/\sqrt2x^{1/4}$);
second, it tells us how to predict approximately when $B(x)$ assumes a
given integer value.

\begin{lemma}
There is a positive constant $C$ such that:
\renewcommand{\labelenumi}{{\it(\alph{enumi})}}
\begin{enumerate}
\item for all real numbers $x\ge1$ and $0\le y\le\min\{x/2,3\sqrt
x\}$, we have
\begin{equation}
\big| (B(x+y)-B(x)) - \frac y{\sqrt2x^{1/4}} \big| < C;  \label{bxbxy}
\end{equation}
\item if we define $z_j = \frac12(3j)^{2/3} - \frac14(3j)^{1/3}$ for
any positive integer $j$, then for all $j>C$ we have $z_j>2$ and
$B((z_j-1)^2) < j < B(z_j^2)$.
\end{enumerate}
\label{Bxlem}
\end{lemma}

If the proof of Lemma~\ref{floorlem} was omitted due to its
tediousness, the proof of this lemma should be omitted and then
buried\dots. The idea of the proof is to rewrite $B_0(x)x^{-3/4}$
using the new variable $t=x^{-1/4}$, and then expand in a Taylor
series in $t$ (a slight but easily overcome difficulty being that the
term $\gamma(x)(1-\gamma(x))$ is not differentiable when
$\sqrt2x^{1/4}$ is an integer). For the proof of part (b), we also
need to rewrite $z_jj^{-2/3}$ using the new variable $u=j^{-1/3}$ and
expand in a Taylor series in $u$. We remark that the constant $C$ in
Lemma~\ref{Bxlem} can be taken to be quite small---in fact, $C=5$ will
suffice.

With these last lemmas in hand, we can dispatch
Corollaries~\ref{asymcor} and~\ref{polycor} in quick succession.

\medskip
\begin{pflike}{Proof of Corollary \ref{asymcor}:}
Let $x>1$ be a real number and define $R(x) = A(x) - 2\sqrt2x^{3/4}/3
- \sqrt x/2$, as in the statement of the corollary. We will describe
how to prove the following more precise statement:
\begin{equation}
R(x) = \big( {2\sqrt2\over3} + g(\sqrt2x^{1/4}) - h(2\sqrt x) \big)
x^{1/4} + R_1(x),
\label{moreprecise}
\end{equation}
where
\begin{equation*}
g(t) = {\{t\}(1-\{t\})\over\sqrt2} \quad\hbox{and}\quad h(t) =
\begin{cases}
\displaystyle {\{t\}\over\sqrt2}, &\rmif
\displaystyle 0\le\{t\}\le\frac12, \\
\displaystyle \sqrt{1-\{t\}} - {1-\{t\}\over\sqrt2}, &\rmif
\displaystyle \frac12\le\{t\}\le1
\end{cases}
\end{equation*}
and $|R_1(x)|<D$ for some constant $D$. The functions $g$ and $h$ are
continuous and periodic with period 1, and are the causes of the
oscillations in the error term $R(x)$. The expression
$g(\sqrt2x^{1/4})$ goes through a complete cycle when $x$ increases by
about $2\sqrt2x^{3/4}$ (one can show this using Taylor expansions yet
again!), which causes the large-scale bounces in the normalized error
term $R(x)x^{-1/4}$ shown in Figure~\ref{rxbsfig} below. Similarly,
the expression $h(2\sqrt x)$ goes through a complete cycle when $x$
increases by about $\sqrt x$, which causes the smaller-scale stutters
shown in the (horizontally magnified) right-hand graph in
Figure~\ref{rxbsfig}.

\begin{figure}[htbf]
\smaller\smaller
\hskip.4cm\psfig{figure=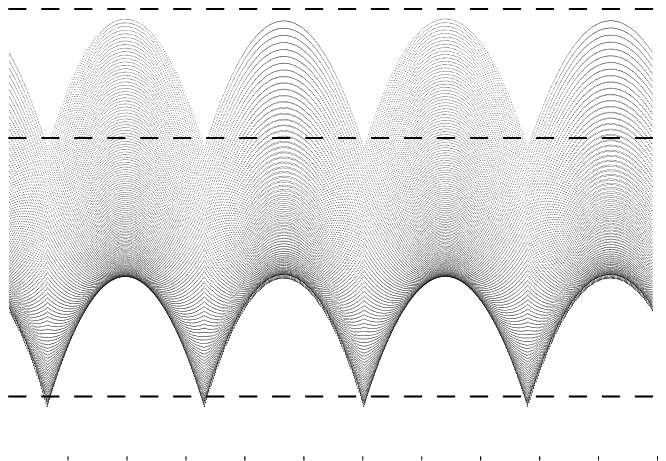}
\put{$x$}{-.2cm}{.15cm}
\put{$R(x)x^{-1/4}$}{-6.8cm}{4.8cm}
{\smaller\smaller
 \put{$10^8$}{-0.75cm}{-.2cm}
 \put{$9.95\cdot10^7$}{-4.17cm}{-.2cm}
 \put{$9.9\cdot10^7$}{-7.03cm}{-.2cm}
 \put{$5/6\sqrt2$}{-8.05cm}{0.6cm}
 \put{$2\sqrt2/3$}{-8.05cm}{3.2cm}
 \put{$19/12\sqrt2$}{-8.35cm}{4.5cm}
}
\hskip.5cm\psfig{figure=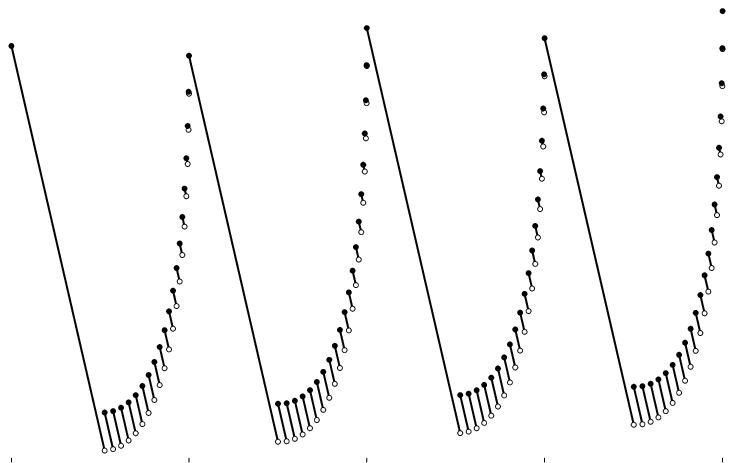}
\put{$x$}{-.1cm}{.2cm}
\put{$R(x)x^{-1/4}$}{-6.8cm}{4.8cm}
{\smaller\smaller
 \put{$800^2$}{-7.8cm}{-.2cm}
 \put{$800\cdot801$}{-6.25cm}{-.2cm}
 \put{$801^2$}{-4.2cm}{-.2cm}
 \put{$801\cdot802$}{-2.65cm}{-.2cm}
 \put{$802^2$}{-0.6cm}{-.2cm}
}
\caption{Big bounces and small stutters for $R(x)x^{-1/4}$}
\label{rxbsfig}
\end{figure}

To establish the formula~(\ref{moreprecise}), we shrewdly add
$B(x)-B(x)$ to the expression defining $R(x)$, which yields
\begin{equation*}
R(x) = (A(x)-B(x)) + \big( \frac{2\sqrt2}3 + {\gamma(x)(1-\gamma(x))
\over \sqrt2} \big) x^{1/4} + B_1(x)
\end{equation*}
from the definition~(\ref{b0def}) of $B_0(x)$. Now $B_1(x)$ is a
bounded function; and since $\gamma(x)=\{\sqrt2x^{1/4}\}$, the
expression $\gamma(x)(1-\gamma(x))/\sqrt2$ is precisely
$g(\sqrt2x^{1/4})$. So what we need to show is that
$B(x)-A(x)=h(2\sqrt x)x^{1/4} + R_2(x)$, where $R_2(x)$ is another
bounded function.

While we won't give all the details, the outline of showing this last
fact is as follows: suppose first that $x\ge m^2$ but that $x$ is less
than the first \asq\ $(m+1+a_{m+1})(m-a_{m+1})$ in the $(2m+1)$-st
flock, so that $A(x)=A(m^2)$. Since $A(m^2)=B(m^2)$ by
Lemma~\ref{am2bm2}, we only need to show that $B(x)-B(m^2)$ is
approximately $h(2\sqrt x)x^{1/4}$; this we can accomplish with the
help of Lemma~\ref{Bxlem}(a).

Similarly, if $x<m^2$ but $x$ is at least as large as the first \asq\
$(m+b_m)(m-b_m)$ in the $2m$-th flock, the same method works as long
as we take into account the difference between $A(m^2)$ and $A(x)$,
which is $\floor{\sqrt{m^2-x}}$ by the Main Theorem. And if $x$ is
close to an \asq\ of the form $m(m-1)$ rather than $m^2$, the same
method applies; even though $A(m(m-1))$ and $B(m(m-1))$ are not
exactly equal, they differ by a bounded amount.

Notice that the functions $g(t)$ and $h(t)$ take values in
$[0,1/4\sqrt2]$ and $[0,1/2\sqrt2]$, respectively. From this and the
formula~(\ref{moreprecise}) we can conclude that
\begin{equation*}
\mathop{\hbox{lim inf}}_{x\to\infty} {R(x)\over x^{1/4}} =
{5\over6\sqrt2} \and \mathop{\hbox{lim sup}}_{x\to\infty} {R(x)\over
x^{1/4}} = {19\over12\sqrt2}.
\end{equation*}
The interested reader can check, for example, that the sequences
$y_j=4j^4+j^2$ and $z_j=(2j^2+j)^2$ satisfy $\lim_{j\to\infty}
R(y_j)/y_j^{1/4} = 5/6\sqrt2$ and $\lim_{j\to\infty}
R(z_j)/z_j^{1/4} = 19/12\sqrt2$.
\end{pflike}

\begin{pflike}{Proof of Corollary \ref{polycor}:}
The algorithms we describe will involve only the following types of
operations: performing ordinary arithmetical calculations $+$, $-$,
$\tms$, $\div$; computing the greatest-integer $\floor\cdot$,
least-integer $\ceil\cdot$, and fractional-part $\{\cdot\}$ functions;
taking square roots; and comparing two numbers to see which is
bigger. All of these operations can be easily performed in polynomial
time. To get the ball rolling, we remark that the functions $a_m$,
$b_m$, and $B(x)$ can all be computed in polynomial time, since their
definitions only involve the types of operations just stated.

We first describe a polynomial-time algorithm for computing the number
of \asqs\ up to a given positive integer $N$. Let $M=\ceil{\sqrt N}$,
so that $(M-1)^2<N\le M^2$. Lemma~\ref{am2bm2} tells us that the
number of \asqs\ up to $M^2$ is $B(M^2)$, and so we simply need to
subtract from this the number of \asqs\ larger than $N$ but not
exceeding $M^2$. This is easy to do by the characterization of \asqs\
given in the Main Theorem. If $N>M(M-1)$, then we want to find the
positive integer $b$ such that $M^2-b^2\le N<M^2-(b-1)^2$, except that
we want $b=0$ if $N=M^2$. In other words, we set
$b=\ceil{M^2-N}$. Then, if $b\le b_M$, the number of \asqs\ up to $N$
is $B(M^2)-b$, while if $b>b_M$, the number of \asqs\ up to $N$ is
$B(M^2)-b_M-1$.

In the other case, where $N\le M(M-1)$, we want to find the positive
integer $a$ such that $(M-a-1)(M+a)\le N<(M-a)(M+a-1)$, except that we
want $a=0$ if $N=M(M-1)$, In other words, we set
$a=\ceil{\sqrt{(M-1/2)^2-N}+1/2}$. Then, if $a\le a_M$, the number of
\asqs\ up to $N$ is $B(M^2)-b_m-1-a$, while if $a>a_M$, the number of
\asqs\ up to $N$ is $B((M-1)^2)$. This shows that $A(N)$ can be
computed in polynomial time, which establishes part (c) of the corollary.

Suppose now that we want to compute the $N$th \asq. We compute in any
way we like the first $C$ \asqs, where $C$ is as in Lemma \ref{Bxlem};
this only takes a constant amount of time (it doesn't change as $N$
grows) which certainly qualifies as polynomial time. If $N\le C$ then
we are done, so assume that $N>C$. Let $M=\ceil{z_N}$, where $z_N$ is
defined as in Lemma \ref{Bxlem}(b), so that $M$ is at least 3 by the
definition of $C$. By Lemma \ref{am2bm2},
\begin{equation*}
A(M^2) = B(M^2) \ge B(z_N^2) > N \quad\hbox{and}\quad A((M-2)^2) =
B((M-2)^2) < B((z_N-1)^2) < N,
\end{equation*}
where the last inequality in each case follows from Lemma
\ref{Bxlem}(b). Therefore the $N$th \asq\ lies between $(M-2)^2$ and
$M^2$, and so is either in the $2M$-th flock or one of the preceding
three flocks. If $0 \le B(M^2)-N \le b_M$, then the $N$th \asq\ is in
the $2M$-th flock, and by setting $b=B(M^2)-N$ we conclude that the
$N$th \asq\ is $M^2-b^2$ and the dimensions of the optimal rectangle
are $(M-b)\tms(M+b)$. If $1+b_M \le B(M^2)-N \le b_M+1+a_M$, then the
$N$th \asq\ is in the $(2M-1)$-st flock, and so on. This establishes
part (d) of the corollary.

Finally, we can determine the greatest \asq\ not exceeding $N$ by
computing $J=A(N)$ and then computing the $J$th \asq, both of which
can be done in polynomial time by parts (c) and (d); and we can
determine whether $N$ is an \asq\ simply by checking whether this
result equals $N$. This establishes the corollary in its entirety.
\end{pflike}

\section{Final Filibuster}\noindent
We have toured some very pretty and precise properties of the \asqs,
and there are surely other natural questions that can be asked about
them, some of which have already been noted. When Grantham posed this
problem, he recalled the common variation on the original
calculus problem where the fence for one of the sides of the rectangle
is more expensive for some reason (that side borders a road or
something), and suggested the more general problem of finding the most
cost-effective rectangle with integer side lengths and area at most
$N$, where one of the sides must be fenced at a higher cost. This
corresponds to replacing $s(n)$ with the more general function
$s_\alpha(n) = \min_{d\mid n}(d+\alpha n/d)$, where $\alpha$ is some
constant bigger than 1. While the elegance of the characterization of
such ``$\alpha$-\asqs'' might not match that of Corollary
\ref{tricor}, it seems reasonable to hope that an enumeration every
bit as precise as the Main Theorem would be possible to establish.

How about generalizing this problem to higher dimensions? For example,
given a positive integer $N$, find the dimensions of the rectangular
box with integer side lengths and volume at most $N$ whose
volume-to-surface area ratio is largest among all such boxes. (It
seems a little more natural to consider surface area rather than the
sum of the box's length, width, and height, but who knows which
problem has a more elegant solution?) Perhaps these ``almost-cubes''
have an attractive characterization analogous to Corollary
\ref{tricor}; almost certainly a result like the Main Theorem, listing
the almost-cubes in order, would be very complicated. And of course
there is no reason to stop at dimension 3.

In another direction, intuitively it seems that numbers with many
divisors are more likely to be \asqs, and the author thought to test
this theory with integers of the form $n!$. However, computations
reveal that the only values of $n\le500$ for which $n!$ is an \asq\
are $n=1,2,3,4,5,6,7,8,10,11,13,15$. Is it the case that these are
the only factorial \asqs? This seems like quite a hard question to
resolve. Perhaps a better intuition about the \asqs\ is that only
those numbers that lie at the right distance from a number of the form
$m^2$ or $m(m-1)$ are \asqs---more an issue of good fortune than of
having enough divisors.

The reader is welcome to contact the author for the Mathematica code
used to calculate the functions related to \asqs\ described in this
paper. With this code, for instance, one can verify that with
8,675,309 square feet of land at his disposal, it is most
cost-effective for Farmer Ted to build a 2,$919'\tms2$,$972'$
supercoop \dots speaking of which, we almost forgot to finish the Farmer
Ted story:

\story{After learning the ways of the \asqs, Farmer Ted went back to
Builders Square, where the salesman viewed the arrival of his
$\R$-rival with trepidation. But Farmer Ted reassured him, ``Don't
worry---I no longer think it's inane to measure fences in $\N$.'' From
that day onward, the two developed a flourishing business
relationship, as Farmer Ted became an integral customer of the store.}

\noindent And that, according to this paper, is that.

\story{{\it Acknowledgements.\/} The author would like to thank
Andrew Granville and the anonymous referees for their valuable
comments which improved the presentation of this paper. The author
would also like to acknowledge the support of National Science
Foundation grant number DMS 9304580 \dots the NSF may or may not wish
to acknowledge this paper.}


\begin{thebibliography}{1}

\bibitem{ConGuy}
J.H.~Conway and R.K.~Guy, {\it The Book of Numbers}, Copernicus, New
York, 1996

\bibitem{dewdney}
A.K.~Dewdney, {\it The (new) Turing Omnibus}, Computer Science Press,
New York, 1993

\bibitem{sqrtn}
S.W.~Golomb, Problem E2491, \emph{Amer. Math. Monthly}~82
(1975), 854--855

\bibitem{pom}
C.~Pomerance, A tale of two sieves, {\it Notices Amer. Math. Soc.}~43
(1996), no.~12, 1473--1485

\bibitem{slo}
N.J.A.~Sloane, An on-line version of the encyclopedia of integer
sequences, {\it Electron. J. Combin.}~1 (1994), feature 1
(electronic)

\bibitem{sloplo}
N.J.A.~Sloane and S.~Plouffe, {\it The Encyclopedia of Integer
Sequences}, Academic Press Inc., San Diego, 1995

\end{thebibliography}
\end{document}